\documentclass{amsart}

\usepackage{amsmath,amssymb,amsthm}
\usepackage{graphicx,tikz}
\newtheorem{theorem}{Theorem}

\newtheorem*{corollary}{Corollary}

\newtheorem*{lemma}{Lemma}

\theoremstyle{definition}

\theoremstyle{remark}

\begin{document}

\title[]{A Wasserstein Inequality and Minimal\\ Green Energy on compact manifolds}
\keywords{Green Energy, Green's function, Coulomb Gas, Wasserstein distance.}
\subjclass[2010]{31B10, 35K05, 49Q20.} 
\thanks{S.S. is supported by the NSF (DMS-1763179) and the Alfred P. Sloan Foundation.}

\author[]{Stefan Steinerberger}
\address{Department of Mathematics, Yale University, New Haven, CT 06511, USA}
\email{stefan.steinerberger@yale.edu}

\begin{abstract}
Let $M$ be a smooth, compact $d-$dimensional manifold, $d \geq 3,$ without boundary and let $G: M \times M \rightarrow \mathbb{R} \cup \left\{\infty\right\}$
denote the Green's function of the Laplacian $-\Delta$ (normalized to have mean value 0). We prove a bound on the cost of transporting
Dirac measures in $\left\{x_1, \dots, x_n\right\} \subset M$ to the normalized volume measure $dx$ in terms of the Green's function of the Laplacian
$$ W_2\left( \frac{1}{n} \sum_{k=1}^{n}{\delta_{x_k}}, dx\right) \lesssim_M  \frac{1}{n^{1/d}} +  \frac{1}{n} \left| \sum_{k, \ell=1 \atop k \neq \ell}^{n}G(x_k, x_{\ell})\right|^{1/2}.$$
We obtain the same result for the Coulomb kernel $G(x,y) = 1/\|x-y\|^{d-2}$ on the sphere $\mathbb{S}^d$, for $d \geq 3$, where we show that
$$ W_2\left(\frac{1}{n} \sum_{k=1}^{n}{ \delta_{x_k}}, dx\right) \lesssim \frac{1}{n^{1/d}} + \frac{1}{n} \left|  \sum_{k, \ell=1 \atop k \neq \ell}^{n}{\left(\frac{1}{\|x_k - x_{\ell}\|^{d-2}} - c_d \right)} \right|^{\frac{1}{2}},$$
where $c_d$ is the constant that normalizes the Coulomb kernel to have mean value 0. We use this to show that minimizers of the discrete Green energy on compact manifolds have optimal rate of convergence $W_2\left( \frac{1}{n} \sum_{k=1}^{n}{\delta_{x_k}}, dx\right)  \lesssim n^{-1/d}$. The second inequality implies the same result for minimizers of the Coulomb energy on $\mathbb{S}^d$ which was recently proven by Marzo \& Mas.
\end{abstract}

\maketitle

\section{Introduction} The problem of distributing points evenly over a compact domain is classical. It is usually ascribed to a 1904 paper of Thomson \cite{thomson} who was concerned with the position of $n$ electrons on $\mathbb{S}^2$ in such a way that
$$ \sum_{i,j=1 \atop i \neq j}^{n}{\frac{1}{\|x_i - x_j\|}} \qquad \mbox{is minimized.}$$
This problem is known to be very hard, the case $n=5$ was only solved very recently \cite{schwartz}. Problems of this flavor have appeared in countless settings, we refer to \cite{beck, brauchart, cohn, dahlberg, fekete, lubotzky, marzo, shub}, the surveys \cite{brauchart2, hardin} and the textbooks \cite{boro, saff}.\\
Recently, Beltr\'an, Corral and Criado del Rey \cite{carlos} introduced the notion of Green energy: if $M$ is a smooth, compact $d-$dimensional manifold without boundary and $G: M \times M \rightarrow \mathbb{R} \cup \left\{\infty\right\}$
is the Green's function of the Laplacian (normalized to have mean value 0), then we can consider the problem of minimizing
$$ \sum_{i,j=1 \atop i \neq j}^{n}{ G(x_i, x_j)} \qquad \mbox{over all sets of}~n~\mbox{points}.$$
This is an exceedingly natural functional insofar as the Green's function is an intrinsic object on the manifold defined as the kernel solving the equation $-\Delta u = f$, i.e.
$$ - \Delta_x \int_{M}{G(x,y) f(y) dy} = f(x).$$
We note that this object has a singularity on the diagonal, we have $G(x,y) \sim 1/|x-y|^{d-2}$ for $d \geq 3$ (and a logarithmic singularity for $d=2$ dimensions). This ensures that having lots of points very close to each other will result in a large Green energy and motivates minimizing it.
The special case of $M=\mathbb{S}^2$ has attracted a lot of attention \cite{bet, fekete, hardin, tsu}, the functional is 
$$ \sum_{i,j=1 \atop i \neq j}^{n}{ \log{\|x_i - x_j\|^{-1}}} \qquad \mbox{and known as the logarithmic energy}.$$
Beltr\'an \cite{beltran} discusses this quantity in terms of a facility location problem.
Beltr\'an, Corral and Criado del Rey \cite{carlos} established that minimizing sets $\left\{x_1, \dots, x_n \right\}$ converge weakly to the Lebesgue measure in the sense that
$$ \lim_{n \rightarrow \infty} \frac{1}{n} \sum_{k=1}^{n}{\delta_{x_k}} \rightharpoonup dx.$$
More recently, Criado del Rey \cite{criado} proved that any minimizing configuration of the Green energy has optimal separation: for all $i \neq j$
$$ |x_i - x_j| \geq \frac{c}{n^{1/d}}.$$
On $\mathbb{R}^n$, the Green's function is merely the Coulomb energy $G(x,y) = c_d/\|x-y\|^{d-2}$ and such sets of points, possibly confined by a growing potential, have been studied intensively,  we refer to \cite{blanc, ch, garcia, leb, rou} and references therein.

\section{Main Result} 
\subsection{Main Result.} The purpose of our paper is to give an inequality connecting this type of quantity to the Wasserstein distance.
The Wasserstein distance \cite{wasser, villani} is a natural notion of distance between measures (roughly `how much mass has to be transported how far from an initial measure to achieve a target measure'). 
We introduce the $p-$Wasserstein distance between two measures
$\mu$ and $\nu$ as
$$ W_p(\mu, \nu) = \left( \inf_{\gamma \in \Gamma(\mu, \nu)} \int_{M \times M}{ |x-y|^p d \gamma(x,y)}\right)^{1/p},$$
where $| \cdot |$ is the usual distance on the torus and $\Gamma(\mu, \nu)$ denotes the collection of all measures on $M \times M$
with marginals $\mu$ and $\nu$, respectively (also called the set of all couplings of $\mu$ and $\nu$). 
Our two measures under consideration are 
$$ \mu = \frac{1}{n} \sum_{k=1}^{n}{\delta_{x_k}} \qquad \mbox{and} \qquad \nu = dx,$$
where $dx$ refers to the normalized volume measure. It is relatively easy to see that, we have an
(optimal) lower bound that is independent of the set $\left\{x_1, \dots, x_n\right\}$
$$W_2\left( \frac{1}{n} \sum_{k=1}^{n}{\delta_{x_k}}, dx\right) \gtrsim \frac{1}{n^{1/d}}.$$
It is now natural to ask how quickly this Wasserstein distance $W_2\left( \frac{1}{n} \sum_{k=1}^{n}{\delta_{x_k}}, dx\right)$ tends to 0 in $n$. It is not difficult to see that if the point set $\left\{x_1, \dots, x_n\right\}$ is not asymptotically uniformly distributed, then the Wasserstein distance does not tend to 0. A quantitative rate of convergence can thus be understood as an improvement over uniform distribution (which is established for many problems of the flavor discussed above and, in particular, for minimizers of the Green energy in \cite{carlos}).\\

Our main result is a general inequality connecting the Wasserstein distance between a discrete set of points and the Green's function. This may be understood as yet manifestation of the phenomenon that Wasserstein distance, like the Green function, operates approximately at Sobolev scale $\dot H^{-1}$ (see \cite{peyre, villani}).
\begin{theorem} Let $M$ be a smooth, compact $d-$dimensional manifold without boundary, $d \geq 3$, and let $G: M \times M \rightarrow \mathbb{R} \cup \left\{\infty\right\}$ denote the Green's function of the Laplacian normalized to have average value 0 over the manifold. Then, for any set of $n$ points
$\left\{x_1, \dots, x_n\right\} \subset M$, we have
$$ W_2\left( \frac{1}{n} \sum_{k=1}^{n}{\delta_{x_k}}, dx\right) \lesssim_M  \frac{1}{n^{1/d}} +  \frac{1}{n} \left| \sum_{k \neq \ell} G(x_k, x_{\ell})\right|^{1/2}.$$
If the manifold is two-dimensional, $d=2$, then we have
$$ W_2\left( \frac{1}{n} \sum_{k=1}^{n}{\delta_{x_k}}, dx\right) \lesssim_M  \frac{\sqrt{\log{n}}}{n^{1/2}} +  \frac{1}{n} \left| \sum_{k \neq \ell} G(x_k, x_{\ell})\right|^{1/2}.$$
\end{theorem}
The result is sharp for $d \geq 3$ and sharp up to possibly the factor of $\sqrt{\log{n}}$ in $d=2$. One way to see this is by computing asymptotics on the Green energy which follows as a byproduct from our approach (this Corollary can be interpreted as related to the work of Wagner \cite{wagner} for Coulomb energy on the sphere).
\begin{corollary} Let $M$ be a smooth, compact $d-$dimensional manifold without boundary and $d \geq 3$, then
$$\sum_{k, \ell =1 \atop k \neq \ell}^{n} G(x_k, x_{\ell}) \gtrsim_M - n^{2-2/d}.$$
If the manifold is two-dimensional, then
$$\sum_{k, \ell =1 \atop k \neq \ell}^{n} G(x_k, x_{\ell}) \gtrsim_M - n^{} \log{n}.$$
\end{corollary}
The Theorem and the Corollary combined then show that for points minimizing the Green energy for $d \geq 3$ we have
$$ W_2\left( \frac{1}{n} \sum_{k=1}^{n}{\delta_{x_k}}, dx\right) \lesssim_M  \frac{1}{n^{1/d}} \qquad \mbox{which is best possible.}$$

This also refines \cite{carlos} showing that minimizers of the Green energy equidistribute on the manifold: if we have an infinite sequence of points for which
$$ \lim_{n \rightarrow \infty} \frac{1}{n^2} \sum_{k, \ell =1 \atop k \neq \ell}^{n} G(x_k, x_{\ell}) = 0,$$
then the sequence is asymptotically uniformly distributed (because the Wasserstein distance tends to 0).\\

This result should have an analogue for the Coulomb energy $1/|x|^{d-2}$ which has Fourier decay $1/|\xi|^2$ corresponding to the space $\dot H^{-1}$ (which in turn is connected to the Wasserstein distance). We make this precise for the sphere $\mathbb{S}^d$. 

\begin{theorem} Let $d \geq 3$ and let $\left\{x_1, \dots, x_n \right\} \subset \mathbb{S}^d$, then
\begin{align*}
 W_2\left(\frac{1}{n} \sum_{k=1}^{n}{ \delta_{x_k}}, dx\right) \lesssim \frac{1}{n^{1/d}} + \left| \frac{1}{n^2} \sum_{k, \ell=1 \atop k \neq \ell}^{n}{\left(\frac{1}{\|x_k - x_{\ell}\|^{d-2}} -c_d \right)}\right|^{\frac{1}{2}},
 \end{align*}
 where 
 $$ c_d =  \frac{1}{|\mathbb{S}^d|^2} \int_{\mathbb{S}^d \times \mathbb{S}^d} \frac{dx dy}{\|x-y\|^{d-2}}$$
 is the constant that normalizes the Coulomb energy to mean 0.
\end{theorem}
If $\left\{x_1, \dots, x_n \right\} \subset \mathbb{S}^d$ is a set of points minimizing the Coulomb energy, then a result of Wagner \cite{wagner} implies that the expression in Theorem 2 is smaller than $ \lesssim c_d n^{-1/d}$ and we obtain the desired result stating that minimizers of the Coulomb energy have optimal convergence rate in $W_2$ -- this optimal convergence speed was very recently proven by Marzo \& Mas \cite[Theorem 1.5]{marzo}; they also ask in \cite{marzo} whether such results are possible for Green energies on manifolds, this is answered by our Theorem 1.
Our proof of Theorem 2 is rather straightforward adaption of the proof of Theorem 1 and uses an equivalent definition of the Sobolev norm $\dot H^{-1}$ in terms of spherical harmonics on the sphere. The special function of the Laplacian of the Riesz kernel is used to compensate for the lack of identities satisfied by the Green function. It is not clear to us whether $W_2$ is the endpoint or whether the optimal rate of convergence can also be proven in $W_{2+\varepsilon}$ for some $\varepsilon>0$.

\subsection{Related results.} Our proof of Theorem 1 follows an underlying philosophy going back to Rougery \& Serfaty \cite{rou}, Chafa$\ddot{\i}$, Hardy \& Ma$\ddot{\i}$da \cite{ch} and Garc\'ia-Zelada \cite{garcia}: we smooth the measure a bit, obtain a regular enough object that allows us to go from the Wasserstein distance to $\dot H^{-1}$ which is the natural scale for the Green energy, and control the errors. Superharmonicity of the Green function plays a role and this motivates the use of the heat kernel on compact manifolds \cite{garcia}. In particular, several of the error estimates we will use can also be found in \cite{garcia}. One notable difference in our use of the heat kernel is an application of Aronson's bound simplifying some of the steps. 
Another notable difference is that the previous results are all phrased in $W_1$ whereas our results are phrased in $W_2$ (using an inequality of Peyr\'e \cite{peyre}). This suggest the possibility that at least some of the previous results may also hold for $W_2$ (this was posed as a problem in \cite{ch}). However, many of the previous results work with bounds for $W_1(\mu, \nu)$ for two arbitrary measures in terms of the energy and this is unlikely to generalize: our argument uses Peyr\'e's inequality which would then result in bounds in terms of the weighted space $\dot H^{-1}(\nu)$. We always work with $\nu = dx$ allowing us to recover $\dot H^{-1}(\nu) = \dot H^{-1}(dx) = \dot H^{-1}$.
 Our proof of Theorem 2 relies on the smoothing procedure in Theorem 1 and then uses, implicitly, the Funk-Hecke formula for spherical harmonics to establish an equivalence between the renormalized Green energy and the Sobolev space $\dot H^{-1}$ (see \cite{marzo}). The Laplacian of the Green's function has, by definition, a particularly simple form and this is no longer true for the Coulomb kernel; however, its Laplacian is simple enough to bootstrap the necessary estimate up to the critical time scale.

\subsection{An Application in Number Theory: Diaphony.} We sketch an interesting application to the case of $d=1$ on the one-dimensional torus $\mathbb{T} = [0,1]$ (with endpoints identified). The Green's function is translation invariant and given by
$$ G(x,y) = G(x-y) = \frac{|x-y|^2}{2} - \frac{|x-y|}{2} + \frac{1}{12}.$$
Given a finite set $\left\{x_1, \dots, x_N\right\}$, we associate the measure 
$$ \mu = \frac{1}{N} \sum_{k=1}^{N}{\delta_{x_k}}.$$
A natural quantity that is frequently studied (see e.g. \cite{dick, drmota, kuipers}) is the discrepancy
$$ D_N(\mu) = \sup_{J \subset \mathbb{T} \atop J~\mbox{\tiny interval}}{ | \mu(J) - |J| |}.$$
It is easy to see that $N^{-1} \leq D_N \leq 1$. The inequality
$$ W_1(\mu, dx) \lesssim D_N(\mu)$$
is easy to see and follows from Monge-Kantorovich duality (see e.g. \cite{stein}). Another notion of regularity is Zinterhof's diaphony \cite{zinterhof} and can be defined as (though it is not usually defined in this manner, see also \cite{drmota} for more details)
\begin{align*}
 F_N(\mu) = \left( \sum_{k \in \mathbb{Z} \atop k \neq 0}{ \frac{| \widehat{\mu}(k)|^2}{k^2}}\right)^{1/2} &= \left( \frac{1}{N^2} \sum_{k, \ell =1}^{N}G(x_k - x_{\ell}) \right)^{1/2} \\
 &= \left( \frac{1}{N^2} \sum_{k, \ell =1}^{N}G(x_k, x_{\ell}) \right)^{1/2}.
 \end{align*}
This establishes a re-interpretation of Zinterhof's diaphony in terms of Green energy and vice versa. Some estimates for the $W_2$ distance of number-theoretic sequence to their equilibrium have been carried out in \cite{stein}. We explicitly note
$$ W_1(\mu, dx) \lesssim D_N(\mu)$$
$$ W_2(\mu, dx) \lesssim F_N(\mu)$$
Existing work on the size of the diaphony of very regularly distributed sets of points that are defined via regular structures in number theory (inversion in finite field, irrational rotations on the torus, ...) immediately carry over as examples for the Green energy.

\section{proofs}
\subsection{Proof of the Theorem.}
 We first give a proof for $d \geq 3$ dimensions to simplify exposition. The necessary changes for $d=2$ dimensions are outlined at the end of the proof. The first step of the argument, control of the Wasserstein under the heat equation is known; the short argument is repeated here for convenience of the reader (the argument would go through for $W_p$ distance for all $1 \leq p < \infty$).
 
 \begin{lemma}[see also \cite{garcia, stein}] Let $\mu$ be a probability measure on the compact manifold $M$. Then
 $$ W_2(\mu, e^{t \Delta} \mu) \lesssim_M \sqrt{t},$$
 where the implicit constant depends only on the manifold.
 \end{lemma}
 \begin{proof} We apply the heat equation for a short time to $\mu$. We interpret the heat equation as convolution with the heat kernel and we interpret the heat kernel as a transport plan.
The result of this
transport plan will be a new mass distribution given by
$$ g(x) = \int_{M}{ e^{t\Delta}\delta_x(y) d\mu(y)} \quad \mbox{at $W_2^2$-cost} \quad \int_{M}{ \int_{M}{ |x-y|^2 e^{t\Delta}\delta_x(y) d\mu(y)} dx}.$$
This transportation cost is easy to bound: we use a classical bound of Aronson \cite{aronson}, that will also play a role
in the proof of our main result,
$$ e^{t\Delta}\delta_x(y) \leq \frac{c_1}{t^{d/2}} \exp \left( -\frac{|x-y|^2}{c_2 t} \right)$$
for some $c_1, c_2$ depending only on $(M,g)$ and obtain
$$ \int_{M}{ \int_{M}{ |x-y|^2 e^{t\Delta}\delta_x(y)d\mu(y)} dx} \lesssim_{(M,g)}  \int_{M}{ \int_{M}{\frac{ |x-y|^2}{t^{d/2}} \exp \left( -\frac{|x-y|^2}{c_2 t} \right) d\mu(y)} dx}.$$
However, it is easily seen that for some universal constants depending on the manifold
$$  \int_{M}{\frac{ |x-y|^2}{t^{d/2}} \exp \left( -\frac{|x-y|^2}{c_2 t} \right) dx} \lesssim_{c_2,(M,g)} t^{}.$$
Altogether, this implies the squared $W_2$ transportation cost is bounded by (since $c_2 \lesssim_{(M,g)} 1$)
$$  \int_{M}{ \int_{M}{ |x-y|^2 e^{t\Delta}\delta_x(y) d\mu(y)} dx} \lesssim_{(M,g)}   t \|\mu\|_{L^1} = t.$$

 \end{proof}

 \begin{proof}[Proof of Theorem 1]
We abbreviate 
$$ \mu = \frac{1}{N} \sum_{k=1}^{N}{ \delta_{x_k}}$$
and first use the triangle inequality
\begin{align*}
 W_2( \mu, dx) &\leq W_2(\mu, e^{t \Delta} \mu) + W_2(e^{t \Delta} \mu, dx),
 \end{align*}
 where $e^{t\Delta}$ is the solution of the heat equation at time $t$. We can define this operator spectrally via
 $$ e^{t\Delta}f = \sum_{k = 0}^{\infty}{ e^{-\lambda_k t} \left\langle f, \phi_k \right\rangle \phi_k},$$
 where $\phi_k$ denotes the sequence of Laplacian eigenfunctions on the manifold, i.e. $-\Delta \phi_k = \lambda_k \phi_k$ (normalized to have $\|\phi_k\|_{L^2} = 1$).
We have already seen in the Lemma above that $ W_2(\mu, e^{t \Delta} \mu) \lesssim \sqrt{t}.$
 For the second term, we use an estimate of Peyr\'e \cite{peyre} that implies
 $$ W_2(e^{t \Delta} \mu, dx) \lesssim \left\| -\frac{1}{\mbox{vol}(M)} + e^{t \Delta} \mu \right\|_{\dot H^{-1}} = \left\|  e^{t \Delta} \mu \right\|_{\dot H^{-1}} ,$$
 where $\dot H^{-1}$ is a Sobolev space whose norm can be defined spectrally via
 $$ \|f\|_{\dot H^{-1}}^2 = \sum_{k=1}^{\infty}{\frac{ \left\langle f, \phi_k \right\rangle^2}{\lambda_k}}.$$
 The one-dimensional case of Peyr\'e's inequality is an identity and can also be found in \cite[Exercise 64]{santa}.
 We note that the Green function is defined in a similar way, i.e. spectrally via
 $$ \int_{M} G(x,y) f(y) dy = \sum_{k=1}^{\infty}{ \frac{ \left\langle f, \phi_k \right\rangle}{\lambda_k}\phi_k(x)}.$$
As a consequence, we have that
\begin{align*}
 \int_{M \times M} G(x,y) f(x) f(y) dx dy &= \left\langle \int_{M} G(x,y) f(y) dy, f(x) \right\rangle \\
 &=  \sum_{k=1}^{\infty}{ \frac{ \left\langle f, \phi_k \right\rangle^2}{\lambda_k}} = \|f\|_{\dot H^{-1}}^2.
 \end{align*}
 We note that the heat equation and the Green function are both spectral multipliers and thus,  whenever $s_1 + t_1 = s_2 + t_2$ and all four numbers are positive,
 \begin{align*}
 \int_{M} \int_{M} G(x,y) e^{s_1\Delta} f(x) e^{t_1 \Delta} g(y) dx dy &= \sum_{k=1}^{\infty}{ e^{-s_1\lambda_k} \frac{\left\langle f, \phi_k\right\rangle \left\langle g, \phi_k \right\rangle}{\lambda_k} e^{-t_1\lambda_k }} \\
  &= \sum_{k=1}^{\infty}{ e^{-s_2\lambda_k} \frac{\left\langle f, \phi_k\right\rangle \left\langle g, \phi_k \right\rangle}{\lambda_k} e^{-t_2\lambda_k }} \\
 &=  \int_{M} \int_{M} G(x,y) e^{s_2\Delta} f(x) e^{t_2 \Delta} g(y) dx dy.
 \end{align*}
We can now write
 \begin{align*}
   \left\| e^{t \Delta} \mu \right\|_{\dot H^{-1}}^2 &= \int_{M} \int_{M} G(x,y) e^{t\Delta} \mu(x) e^{t \Delta} \mu (y) dx dy\\
   &= \frac{1}{N^2} \sum_{k, \ell} \int_{M} \int_{M} G(x,y) e^{t\Delta} \delta_{x_k}(x) e^{t \Delta} \delta_{x_{\ell}} (y) dx dy \\
   &= \frac{1}{N^2} \sum_{k} \int_{M} \int_{M} G(x,y) e^{t\Delta} \delta_{x_k}(x) e^{t \Delta} \delta_{x_{k}} (y) dx dy\\
   &+ \frac{1}{N^2} \sum_{k \neq \ell} \int_{M} \int_{M} G(x,y) e^{t\Delta} \delta_{x_k}(x) e^{t \Delta} \delta_{x_{\ell}} (y) dx dy.
   \end{align*}
We use the self-adjointness of spectral multipliers (and both convolution with $G$ as well as the heat kernel are spectral multipliers, moreover the heat kernel is a semigroup and $e^{t\Delta} e^{t\Delta} = e^{2t \Delta}$) and rewrite the first term as
\begin{align*}
 \int_{M} \int_{M} G(x,y) e^{t\Delta} \delta_{x_k}(x) e^{t \Delta} \delta_{x_{k}} (y) dx dy &= \int_{M} \int_{M} G(x,y) \delta_{x_k}(x) e^{2t \Delta} \delta_{x_{k}} (y) dx dy \\
 &= \int_{M} G(x_k,y) e^{2t \Delta} \delta_{x_{k}} (y) dy
 \end{align*}  
 We have a very good understanding of the heat kernel since
 $$ e^{2t \Delta} \delta_{x_{k}} (y) \lesssim \begin{cases} t^{-d/2} \qquad &\mbox{if}~d_g(x,y) \leq \sqrt{t} \\
 \mbox{exponentially decaying} \qquad &\mbox{otherwise.} \end{cases}$$
 More formally, we use a classical bound of Aronson \cite{aronson} for the heat kernel on manifolds stating that for constants $c_1, c_2$ depending on the manifold
$$ e^{t\Delta} \delta_x (y) \leq \frac{c_1}{t^{d/2}} \exp \left( -\frac{|x-y|^2}{c_2 t} \right),$$
where $|x-y| = d_g(x,y)$ denotes the geodesic distance (we will use this notation henceforth).
 We can couple this with the estimate (see e.g. Aubin \cite{aubin})
 $$ G(x,y) \lesssim_M \frac{1}{|x-y|^{d-2}}$$
 and obtain
 \begin{align*}
   \int_{M} G(x_k,y) e^{2t \Delta} \delta_{x_{k}} (y) dy &\lesssim \int_{\mathbb{R}^n}^{}{ \frac{e^{2t \Delta} \delta_{0} (x)}{|x|^{d-2}}dx} \\
   &\lesssim \int_{0}^{\infty}{ \frac{c_1}{t^{d/2}} \frac{ \exp \left( -\frac{r^2}{c_2 t} \right)}{r^{n-2}} r^{n-1}} dr \lesssim t^{1-\frac{d}{2}}.
   \end{align*}
This implies
$$ \frac{1}{N^2} \sum_{k} \int_{M} \int_{M} G(x,y) e^{t\Delta} \delta_{x_k}(x) e^{t \Delta} \delta_{x_{k}} (y) dx dy \lesssim \frac{ t^{1-\frac{d}{2}}}{N}.$$

It remains to bound the second term. We can again use that Fourier multipliers commute to argue that
 \begin{align*}
  \int_{M} \int_{M} G(x,y) e^{t\Delta} \delta_{x_k}(x) e^{t \Delta} \delta_{x_{\ell}} (y) dx dy &=  \int_{M} \int_{M} G(x,y) \delta_{x_k}(x) e^{2t \Delta} \delta_{x_{\ell}} (y) dx dy \\
  &=  \int_{M} G(x_k ,y)e^{2t \Delta} \delta_{x_{\ell}} (y) dy.
 \end{align*}
 We understand this value for $t$ very small since the Green function is integrable and
 $$ \lim_{t \rightarrow 0} \int_{M} G(x_k ,y)e^{2t \Delta} \delta_{x_{\ell}} (y) dy = G(x_k, x_{\ell})$$
 and will now control the variation in time. We note that if $x_k = x_{\ell}$ for some $k \neq \ell$, then our upper bound is infinity/undefined and the entire statement is vacuous. We can thus assume $x_k \neq x_{\ell}$ for $k \neq \ell$. 
 The heat kernel solves the heat equation and thus
 $$ \frac{\partial}{\partial t}  e^{t \Delta} \delta_{x_{\ell}} (y) = \Delta_y e^{t \Delta} \delta_{x_{\ell}} (y)$$
 which we use in combination with
 $$ \Delta_y G(x, y) = \frac{1}{\mbox{vol}(M)} - \delta_{x} \qquad \mbox{and} \qquad \int_{M} e^{2t \Delta} \delta_{x_{\ell}} (y) dy = 1,$$
 where the first identity is in the sense of distributions (and will be used paired against a smooth function).
 We write

 \begin{align*}
 \frac{1}{2} \frac{\partial}{\partial t}  \int_{M} G(x_k ,y)e^{2t \Delta} \delta_{x_{\ell}} (y) dy &=   \int_{M} G(x_k ,y) \Delta_y e^{2t \Delta} \delta_{x_{\ell}} (y) dy\\
 &=   \int_{M} \Delta_y G(x_k ,y)  e^{2t \Delta} \delta_{x_{\ell}} (y) dy \\
 &=  \int_{M} \left(\frac{1}{\mbox{vol}(M)} - \delta_{x_k}\right)  e^{2t \Delta} \delta_{x_{\ell}} (y) dy\\
 &= \frac{1}{\mbox{vol}(M)} - \left( e^{2 t \Delta} \delta_{x_{\ell}} \right)(x_k)\\
 &\leq  \frac{1}{\mbox{vol}(M)} 
 \end{align*}
This implies
 $$ \frac{1}{N^2} \sum_{k \neq \ell} \int_{M} \int_{M} G(x,y) e^{t\Delta} \delta_{x_k}(x) e^{t \Delta} \delta_{x_{\ell}} (y) dx dy \leq \frac{2t}{\mbox{vol}(M)} + \frac{1}{N^2} \sum_{k \neq \ell} G(x_k, x_{\ell}).$$
 Altogether, collecting all the estimates, we have
 \begin{align*}
   W_2( \mu, dx) &\lesssim W_2( \mu, e^{t\Delta} \mu) + W_2(e^{t\Delta}\mu, dx) \\
   &\lesssim \sqrt{t} + \left( \frac{t^{1-\frac{d}{2}}}{N} + \frac{2t}{\mbox{vol}(M)} + \frac{1}{N^2} \sum_{k \neq \ell} G(x_k, x_{\ell}) \right)^{1/2}.
   \end{align*}
 Setting $t = N^{-2/d}$ results in
 $$ W_2\left(\frac{1}{N} \sum_{k=1}^{N}{ \delta_{x_k}}, dx\right) \lesssim_M \frac{1}{N^{1/d}} +\frac{1}{N}  \left|\sum_{k \neq \ell} G(x_k, x_{\ell}) \right|^{1/2}.$$
 We observe that the Green energy may actually be negative and, at first glance, it may look like one could get improved results. However, we will show in the proof of the corollary that 
 $$ \frac{1}{N}  \left|\sum_{k \neq \ell} G(x_k, x_{\ell}) \right|^{1/2} \gtrsim_M - N^{-1/d}$$
 and thus our application of the triangle inequality in this form is not lossy.
 We quickly note the necessary changes for the case $d=2$. We observe that in that case
 $$ G(x,y) \lesssim |\log{|x-y|}|$$
 and thus, for the first term,
  \begin{align*}
   \int_{M} G(x_k,y) e^{2t \Delta} \delta_{x_{k}} (y) dy &\lesssim \int_{\mathbb{R}^n}^{}{ |\log{|x|}|e^{2t \Delta} \delta_{0} (x) dx} \\
   &\lesssim \int_{0}^{\infty}{ \frac{c_1}{t^{}} |\log{r}| \exp \left( -\frac{r^2}{c_2 t} \right) r dr} \lesssim \log{(1/t)}.
   \end{align*}
 The off-diagonal term behaves exactly as before and we obtain
  $$  W_2( \mu, dx) \lesssim_M  \sqrt{t} + \left( \frac{\log{(1/t)}}{N} + t + \frac{1}{N^2} \sum_{k \neq \ell} G(x_k, x_{\ell}) \right)^{1/2}.$$
Setting $t=1/N$ results in the desired statement.
\end{proof}

\subsection{Proof of the Corollary}
The Corollary follows from the argument developed in Theorem 1. We again first deal with the case $d \geq 3$ and then discuss the necessary modifications to deal with $d=2$.

\begin{proof} We make use of the trivial identity
$$    \left\| e^{t \Delta} \mu \right\|_{\dot H^{-1}}^2 \geq 0.$$
At the same time, we can control its expansion in terms of Green's function and the bounds obtained in Theorem 1 via
 \begin{align*}
   0 \leq \left\| e^{t \Delta} \mu \right\|_{\dot H^{-1}}^2 &\leq \frac{1}{N^2} \sum_{k} \int_{M} \int_{M} G(x,y) e^{t\Delta} \delta_{x_k}(x) e^{t \Delta} \delta_{x_{k}} (y) dx dy\\
   &+ \frac{1}{N^2} \sum_{k \neq \ell} \int_{M} \int_{M} G(x,y) e^{t\Delta} \delta_{x_k}(x) e^{t \Delta} \delta_{x_{\ell}} (y) dx dy.
   \end{align*}
   We have the inequalities
   \begin{align*}
   \frac{1}{N^2} \sum_{k} \int_{M} \int_{M} G(x,y) e^{t\Delta} \delta_{x_k}(x) e^{t \Delta} \delta_{x_{k}} (y) dx dy &\lesssim \frac{t^{1-\frac{d}{2}}}{N}\\
  \frac{1}{N^2} \sum_{k \neq \ell}   \int_{M} \int_{M} G(x,y) e^{t\Delta} \delta_{x_k}(x) e^{t \Delta} \delta_{x_{\ell}} (y) dx dy &\leq \frac{t}{\mbox{vol}(M)} + \frac{1}{N^2} \sum_{k \neq \ell} G(x_k, x_{\ell})
   \end{align*}
Setting $t= N^{-2/d}$ results in the desired bound. Using the modified bounds yields in the result for $d=2$.
\end{proof}

We quickly sketch how one would derive, via an explicit construction, for $d \geq 3$ sets of points that establish
$$ \inf_{x_1, \dots, x_n} \sum_{k \neq \ell} G(x_k, x_{\ell}) \lesssim -n^{2-2/d}.$$
 Our construction is asymptotic as $n \rightarrow \infty$ (so that our considerations are local instead of global) and we illustrate it for the flat torus $\mathbb{T}^d$.
We assume that $n$ is of the form $n = m^d$ and partition the manifold into boxes $(Q_j)_{j=1}^{n}$ of sidelength $m^{-1} = n^{-1/d}$ and put points in the middle of these cubes. We fix a point $x_1$ and want to estimate the expression
$ \sum_{k=2}^{n}{G(x_1, x_{k)}}.$
We start by remarking that
$$ 0 = \int_{\mathbb{T}^d}{G(x_1, y) dy} = \int_{Q_1}{ G(x_1, y) dy} + \sum_{j=2}^{n} \int_{Q_j}{ G(x_1, y) dy},$$
where $Q_1$ is the cube containing $x_1$. We observe that
$$ \left|  \int_{Q_1}{ G(x_1, y) dy}  \right| \lesssim \int_0^{n^{-1/d}} \frac{1}{r^{n-2}} r^{n-1} dr = n^{-2/d}.$$
As for the remaining cubes, we can expand 
\begin{align*}
 G(x_1, y) &= G(x_1, x_k) + \nabla G(x_1, y)\big|_{y=x_k} (y-x_1) \\
 &+ \frac{1}{2} \left\langle y-x_k, D^2 G(x_1, y)\big|_{y=x_k}(y-x_k) \right\rangle + \mbox{l.o.t.}
 \end{align*}
Integrating this over $Q_k$ results in cancellation of the linear term and yields
\begin{align*}
 \int_{Q_j}{ G(x_1, y) dy} - \frac{1}{n} G(x_1, x_k) &= \int_{Q_k}  \frac{1}{2} \left\langle y-x_k, D^2 G(x_1, y)\big|_{y=x_k}(y-x_k) \right\rangle\\
 &\sim n^{-2/d} |Q_j|.
 \end{align*}
Multiplying by $n$, summing over all $n-1$ cubes and then adding over all other $n$ points results in the desired scaling. One could see how this, in the general case where a subdivision into boxes may not be feasible, can be interpreted as a facility location problem (see B\'eltran \cite{beltran}).

\subsection{Proof of Theorem 2.}
The new ingredient is a folklore Lemma that relates the Green function on the sphere to the Coulomb kernel by showing that the Coulomb kernel has exactly the appropriate decay in terms of Fourier coefficients. That computation is carried out using the Funk-Hecke formula, a recent and very nice presentation can be found in Marzo \& Mas \cite{marzo}.
\begin{lemma}[Folklore, see e.g. \cite{marzo}] Let $f \in L^2(\mathbb{S}^d)$. Then
\begin{align*} \|f\|_{\dot H^{-1}}^2 &=  \int_{\mathbb{S}^d}  \int_{\mathbb{S}^d}  G(x,y)f(x) f(y) dx dy\\
&\sim_d \int_{\mathbb{S}^d}  \int_{\mathbb{S}^d}  \frac{f(x) f(y)}{\|x-y\|^{d-2}} dx dy - \left(\frac{1}{|\mathbb{S}^d|}\int_{\mathbb{S}^d}{ f(x) dx} \right)^2  \int_{\mathbb{S}^d \times \mathbb{S}^d}{ \frac{dx dy}{\|x-y\|^{d-2}}}.
\end{align*}
\begin{proof} We note that $G$ is normalized to have mean value and thus annihilates constants. Rewriting $f$ as a constant function and a function having mean value 0 results in
$$ f = \frac{1}{|\mathbb{S}^d|}\int_{\mathbb{S}^d}{ f(x) dx}  + \left(f - \frac{1}{|\mathbb{S}^d|}\int_{\mathbb{S}^d}{ f(x) dx} \right)$$
and plugging this into the right-hand side shows that it remains to deal with functions with mean value 0.
It then suffices to test the equivalence on Laplacian eigenfunctions (spherical harmonics) since
$$ f = \sum_{k=1}^{\infty}{a_k \phi_k}$$
results in
$$ \int_{\mathbb{S}^d}  \int_{\mathbb{S}^d}  G(x,y)f(x) f(y) dx dy = \sum_{k=1}^{\infty}{\frac{a_k^2}{\lambda_k}}.$$
However, we also have
$$ \int_{\mathbb{S}^d}  \int_{\mathbb{S}^d}  \frac{f(x) f(y)}{\|x-y\|^{d-2}} dx dy = \sum_{k, \ell=1}^{\infty} a_k a_{\ell} \int_{\mathbb{S}^d}  \int_{\mathbb{S}^d}  \frac{\phi_k(x) \phi_{\ell}(y)}{\|x-y\|^{d-2}} dx dy $$
and the desired result now follows from the orthogonality relation

$$  \int_{\mathbb{S}^d}  \int_{\mathbb{S}^d}  \frac{\phi_k(x) \phi_{\ell}(y)}{\|x-y\|^{d-2}} dx dy = \begin{cases} 0 \qquad &\mbox{if}~k \neq \ell \\
 \sim_d  \lambda_k^{-1} \qquad &\mbox{if}~k=l \end{cases}.$$
This computation is carried out using the Funk-Hecke formula (see e.g. \cite{marzo}).
\end{proof}
\end{lemma}

\begin{proof}[Proof of Theorem 2] We set again
$$ \mu = \frac{1}{N} \sum_{k=1}^{N}{ \delta_{x_k}}$$
and use, as before,
\begin{align*}
 W_2( \mu, dx) &\leq \sqrt{t} + W_2(e^{t \Delta} \mu, dx) \\
 &\lesssim \sqrt{t} + \left\| e^{t \Delta} \mu \right\|_{\dot H^{-1}}.
 \end{align*}
We use the equivalence of the Green energy and the Riesz kernel discussed in the Lemma and obtain (note that the quantity on the right-hand side is always positive)
\begin{align*} \left\| e^{t \Delta} \mu \right\|_{\dot H^{-1}}^2 &\lesssim_d \frac{1}{N^2} \sum_{k, \ell=1}^{N} \int_{\mathbb{S}^d}  \int_{\mathbb{S}^d}  \frac{ e^{t\Delta} \delta_{x_k}(x) e^{t\Delta} \delta_{x_{\ell}}(x)(y)}{\|x-y\|^{d-2}} dx dy\\
&- \frac{1}{|\mathbb{S}^d|^2} \int_{\mathbb{S}^d \times \mathbb{S}^d}{ \frac{dx dy}{\|x-y\|^{d-2}}}.
\end{align*}
We split in diagonal and off-diagonal terms, the same computation as in the proof of Theorem 1 (where the only property of the Green's function that we used in this step was $|G(x,y)| \lesssim 1/\|x-y\|^{d-2}$) implies
$$  \frac{1}{N^2} \sum_{k=1}^{N} \int_{\mathbb{S}^d}  \int_{\mathbb{S}^d}  \frac{ e^{t\Delta} \delta_{x_k}(x) e^{t\Delta} \delta_{x_{k}}(y)}{\|x-y\|^{d-2}} dx dy \lesssim \frac{t^{1-\frac{d}{2}}}{N}.$$
It remains to deal with the off-diagonal terms: since everything can be interpreted as a spectral multiplier, we can write
\begin{align*} \int_{\mathbb{S}^d}  \int_{\mathbb{S}^d}  \frac{ e^{t\Delta} \delta_{x_{\ell}}(x) e^{t\Delta} \delta_{x_{k}}(y)}{\|x-y\|^{d-2}} dx dy &=
\int_{\mathbb{S}^d}  \int_{\mathbb{S}^d}  \frac{  \delta_{x_k}(x) e^{2t\Delta} \delta_{x_{\ell}}(y)}{\|x-y\|^{d-2}} dx dy\\
&=   \int_{\mathbb{S}^d}  \frac{  e^{2t\Delta} \delta_{x_{\ell}}(y)}{\|x_k-y\|^{d-2}} dy.
\end{align*}
We understand this expression as $t=0$ and try to understand its variation in time. As before, we differentiate in time and obtain
\begin{align*}
\frac{1}{2} \frac{\partial}{\partial t} \int_{\mathbb{S}^d}  \frac{  e^{2t\Delta} \delta_{x_{\ell}}(y)}{\|x_k-y\|^{d-2}} dy &=\int_{\mathbb{S}^d}  \frac{1}{{\|x_k-y\|^{d-2}} } \Delta_y e^{2t\Delta} \delta_{x_{\ell}}(y)dy \\
&= \int_{\mathbb{S}^d}  \Delta_y \left( \frac{1}{{\|x_k-y\|^{d-2}} } \right)e^{2t\Delta} \delta_{x_{\ell}}(y)dy.
\end{align*}
In the proof of Theorem 1, we used that the Laplacian of the Green's function has a particularly simple form. This is different here but the Laplacian of a Coulomb kernel is still explicit and we have, for some positive constants $c_1, c_2$ depending only on the dimension (and could be found, for example, in \cite{marzo}),
$$ \Delta_y  \frac{1}{{\|x_k-y\|^{d-2}} } =  c_1 \frac{1}{\|x_k - y\|^{d-2}}  - c_2 \delta_{x_k}.$$
This implies
$$ \frac{1}{2} \frac{\partial}{\partial t} \int_{\mathbb{S}^d}  \frac{  e^{2t\Delta} \delta_{x_{\ell}}(y)}{\|x_k-y\|^{d-2}} dy \leq c_1  \int_{\mathbb{S}^d}  \frac{  e^{2t\Delta} \delta_{x_{\ell}}(y)}{\|x_k-y\|^{d-2}} dy$$
and therefore
$$  \int_{\mathbb{S}^d}  \frac{  e^{2t\Delta} \delta_{x_{\ell}}(y)}{\|x_k-y\|^{d-2}} dy \leq \frac{e^{c_1 t} }{\|x_k - x_{\ell}\|^{d-2}}.$$
This allows us to bound the off-diagonal terms as
$$ \frac{1}{N^2} \sum_{k, \ell=1 \atop k \neq \ell}^{N} \int_{\mathbb{S}^d}  \int_{\mathbb{S}^d}  \frac{ e^{t\Delta} \delta_{x_k}(x) e^{t\Delta} \delta_{x_{\ell}}(x)(y)}{\|x-y\|^{d-2}} dx dy \leq  \frac{e^{c_1 t}}{N^2} \sum_{k, \ell=1 \atop k \neq \ell}^{N} \frac{1}{\|x_k-x_{\ell}\|^{d-2}}.$$
Let us now fix the abbreviation
$$ X=  \sum_{k, \ell=1 \atop k \neq \ell}^{N}{\left(\frac{1}{\|x_k - x_{\ell}\|^{d-2}} -c_d \right)},$$
where $c_d$ is chosen so that the Coulomb kernel has average value zero over the sphere.
Then, for $0 < t < c_1^{-1}$, we use
\begin{align*}
 \frac{1}{N^2} \sum_{k, \ell=1 \atop k \neq \ell}^{N} \int_{\mathbb{S}^d}  \int_{\mathbb{S}^d}  \frac{ e^{t\Delta} \delta_{x_k}(x) e^{t\Delta} \delta_{x_{\ell}}(x)(y)}{\|x-y\|^{d-2}} dx dy &\leq  \frac{e^{c_1 t}}{N^2} \sum_{k, \ell=1 \atop k \neq \ell}^{N} \frac{1}{\|x_k-x_{\ell}\|^{d-2}}\\
 &= \frac{e^{c_1 t}}{N^2} \left( \frac{ N^2}{|\mathbb{S}^d|^2} \int_{\mathbb{S}^d \times \mathbb{S}^d}{ \frac{dx dy}{\|x-y\|^{d-2}}}   + |X|\right)\\
 &\leq \frac{1}{|\mathbb{S}^d|^2} \int_{\mathbb{S}^d \times \mathbb{S}^d}{ \frac{dx dy}{\|x-y\|^{d-2}}} \\
 &+  e c_1 t\frac{ N^2}{|\mathbb{S}^d|^2} \int_{\mathbb{S}^d \times \mathbb{S}^d}{ \frac{dx dy}{\|x-y\|^{d-2}}} \\
 &+ \frac{e^{c_1 t}}{N^2}|X|
\end{align*}
We pick $t = N^{-2/d}$ and thus obtain the estimate
\begin{align*}
\frac{1}{N^2} \sum_{k, \ell=1 \atop k \neq \ell}^{N} \int_{\mathbb{S}^d}  \int_{\mathbb{S}^d}  \frac{ e^{t\Delta} \delta_{x_k}(x) e^{t\Delta} \delta_{x_{\ell}}(x)(y)}{\|x-y\|^{d-2}} dx dy &\leq  \frac{1}{|\mathbb{S}^d|^2} \int_{\mathbb{S}^d \times \mathbb{S}^d}{ \frac{dx dy}{\|x-y\|^{d-2}}}\\
&+ e c_1 \frac{ N^{2-2/d}}{|\mathbb{S}^d|^2} \int_{\mathbb{S}^d \times \mathbb{S}^d}{ \frac{dx dy}{\|x-y\|^{d-2}}}\\
&+ e c_1 \frac{1}{N^2}|X|
\end{align*}
Collecting all the estimates now results in
 $$\left\| e^{t \Delta} \mu \right\|_{\dot H^{-1}}^2 \lesssim N^{-2/d} + \frac{|X|}{N^2}$$
 and thus
 \begin{align*}
 W_2\left( \frac{1}{N} \sum_{k=1}^{N}{ \delta_{x_k}}, dx\right) &\lesssim N^{-1/d} + \frac{1}{N} |X|^{1/2}\\
 &= N^{-1/d} + \frac{1}{N}\left|  \sum_{k, \ell=1 \atop k \neq \ell}^{N}{\left(\frac{1}{\|x_k - x_{\ell}\|^{d-2}} -c_d \right)}\right|^{\frac{1}{2}}
 \end{align*}
 as desired.
\end{proof}

We conclude by remarking that, as in the proof of the Corollary, we can use the argument to show that for any $\left\{x_1, \dots, x_N \right\} \subset \mathbb{S}^d$, for $ d\geq 3$, we have
$$ \sum_{k, \ell=1 \atop k \neq \ell}^{N}{ \frac{1}{\|x_k - x_{\ell}\|^{d-2}}} - \frac{ N^2}{|\mathbb{S}^d|^2} \int_{\mathbb{S}^d \times \mathbb{S}^d}{ \frac{dx dy}{\|x-y\|^{d-2}}} \geq - c_d N^{2-2/d}$$
which is a result of Wagner \cite{wagner}.

\begin{proof}[Sketch of the Argument.] We define, as above,
$$ \mu = \sum_{k=1}^{N}{\delta_{x_k}} \qquad \mbox{and use} \qquad \| e^{t\Delta} \mu \|_{\dot H^{-1}}^2 \geq 0.$$
In short, we use
$$  \| e^{t\Delta} \mu \|_{ H^{-1}}^2  = \mbox{diagonal} + \mbox{off-diagonal},$$
and we have shown above that
$$  \mbox{diagonal} \leq c_1 \frac{t^{1-d/2}}{N} \leq c_1 N^{2-2/d},$$
where $c_1$ depends only on the dimension. We have also seen that, for constants $c_2, c_3$ depending only on the dimension
$$ \mbox{off-diagonal} \leq c_2 e^{c_3 N^{-2/d}} \sum_{k, \ell=1 \atop k \neq \ell}^{N}{ \frac{1}{\|x_k - x_{\ell}\|^{d-2}}}.$$
Suppose now,
$$ \sum_{k, \ell=1 \atop k \neq \ell}^{N}{ \frac{1}{\|x_k - x_{\ell}\|^{d-2}}} = \frac{ N^2}{|\mathbb{S}^d|^2} \int_{\mathbb{S}^d \times \mathbb{S}^d}{ \frac{dx dy}{\|x-y\|^{d-2}}} - c_4 N^{2-2/d},$$
for some $c_4$ sufficiently large (depending on $c_1, c_2, c_3$), then
$$ \| e^{t\Delta} \mu \|_{ H^{-1}}^2 <  \frac{ N^2}{|\mathbb{S}^d|^2} \int_{\mathbb{S}^d \times \mathbb{S}^d}{ \frac{dx dy}{\|x-y\|^{d-2}}}$$
and then
$$ \| e^{t\Delta} \mu \|_{\dot H^{-1}}^2 < 0 \qquad \mbox{which is a contradiction.}$$
\end{proof}

\end{document}